\documentclass[fleqn,preprint,3p,a4paper]{elsarticle}

\usepackage{amssymb}
\usepackage{amsmath}
\usepackage{amsthm}
\usepackage{dcolumn}
\usepackage{endnotes}
\usepackage{tabularx}
\usepackage[matrix,arrow]{xy}
\usepackage{wasysym}

\newtheorem{theorem}{Theorem}[section]
\newtheorem{proposition}[theorem]{Proposition}
\newtheorem{lemma}[theorem]{Lemma}
\newtheorem{corollary}[theorem]{Corollary}

\theoremstyle{definition}
\newtheorem{definition}[theorem]{Definition}
\newtheorem{example}[theorem]{Example}
\newtheorem{remark}[theorem]{Remark}
\newtheorem{question}[theorem]{Question}

\newcommand{\ir}{{\mathsf{Irr}}}

\newcommand{\cl}{{\rm cl}}
\newcommand{\ii}{{\rm int}}

\newcommand{\ua}{\mathord{\uparrow}}
\newcommand{\da}{\mathord{\downarrow}}

\newcommand{\md}{\mathord{\mathsf{D}}}

\newcommand{\kf}{\mathord{\mathsf{RD}}}

\journal{Topology and its Applications}

\begin{document}

\begin{frontmatter}



\title{A direct approach to $K$-reflections of $T_0$ spaces\tnoteref{t1}}
\tnotetext[t1]{This research was supported by the National Natural Science Foundation of China (11661057) and the Natural Science Foundation of Jiangxi Province , China (20192ACBL20045)}

\author[X. Xu]{Xiaoquan Xu}
\ead{xiqxu2002@163.com}
\address[X. Xu]{School of Mathematics and Statistics,
Minnan Normal University, Zhangzhou 363000, China}

\begin{abstract}
  In this paper, we provide a direct approach to $\mathbf{K}$-reflections of $T_0$ spaces. For a full subcategory $\mathbf{K}$ of the category of all $T_0$ spaces and a $T_0$ space $X$, let $\mathbf{K}(X)=\{A\subseteq X : A$ is closed and for any continuous mapping $f : X\longrightarrow Y$
 to a  $\mathbf{K}$-space $Y$, there exists a unique $y_A\in Y$ such that $\overline{f(A)}=\overline{\{y_A\}}\}$ and $P_H(\mathbf{K}(X))$ the space of $\mathbf{K}(X)$ endowed with the lower Vietoris topology.  It is  proved that if $P_H(\mathbf{K}(X))$ is a $\mathbf{K}$-space,  then the pair $\langle X^k=P_H(\mathbf{K}(X)), \eta_X\rangle$, where $\eta_X :X\longrightarrow X^k$, $x\mapsto\overline{\{x\}}$, is the $\mathbf{K}$-reflection of $X$. We call $\mathbf{K}$ an adequate category if for any $T_0$ space $X$, $P_H(\mathbf{K}(X))$ is a $\mathbf{K}$-space. Therefore, if $\mathbf{K}$ is adequate, then $\mathbf{K}$ is reflective in $\mathbf{Top}_0$. It is shown that the category of all sober spaces, that of all $d$-spaces, that of all well-filtered spaces and the Keimel and Lawson's category are all adequate, and hence are all reflective in $\mathbf{Top}_0$. Some major properties of $\mathbf{K}$-spaces and $\mathbf{K}$-reflections of $T_0$ spaces are investigated. In particular, it is proved that if $\mathbf{K}$ is adequate, then the $\mathbf{K}$-reflection preserves finite products of $T_0$ spaces. Our study also leads to a number of problems, whose answering will deepen our understanding of the related spaces and their categorical structures.
\end{abstract}





\end{frontmatter}


\section{Introduction}

Domain theory plays a foundational role in denotational semantics of programming languages. In domain theory, the $d$-spaces, well-filtered spaces and sober spaces form three of the most important classes (see [3-17, 19-29]). Let $\mathbf{Top}_0$ be the category of all $T_0$ spaces and $\mathbf{Sob}$ the category of all sober spaces. Denote the category of all $d$-spaces and that of all well-filtered spaces respectively by $\mathbf{Top}_d$ and $\mathbf{Top}_w$. It is well-known that $\mathbf{Sob}$ is reflective in $\mathbf{Top}_0$ (see \cite{redbook, Jean-2013}). Using $d$-closures, Wyler \cite{Wyler} proved that $\mathbf{Top}_d$ is reflective in $\mathbf{Top}_0$ (see also \cite{Keimel-Lawson, ZhangLi}). Later, Ershov \cite{Ershov_1999} showed that the $d$-completion of $X$ (i.e., the $d$-reflection of $X$) can be obtained by adding the closure of directed sets onto $X$ (and then repeating this process by transfinite induction). In \cite{Keimel-Lawson}, Keimel and Lawson proved that for a full subcategory $\mathbf{K}$ of $\mathbf{Top}_0$ containing $\mathbf{Sob}$, if $\mathbf{K}$ has certain properties, then $\mathbf{K}$ is reflective in $\mathbf{Top}_0$. They showed that $\mathbf{Top}_d$ and some other categories have such properties.   For quite a long time, it is not known whether $\mathbf{Top}_w$ is reflective in $\mathbf{Top}_0$. Recently, following Keimel and Lawson's method, which originated from Wyler's method, Wu, Xi, Xu and Zhao \cite{wu-xi-xu-zhao-19} gave a positive answer to the above problem. Following Ershov's method of constructing the $d$-completion of $T_0$ spaces, Shen, Xi, Xu and Zhao presented a construction of the well-filtered reflection of $T_0$ spaces.

In this paper, we will provide a direct approach to $\mathbf{K}$-reflections of $T_0$ spaces. For a full subcategory $\mathbf{K}$ of $\mathbf{Top}_0$ containing $\mathbf{Sob}$ and a $T_0$ space $X$, let $\mathbf{K}(X)=\{A\subseteq X : A$ is closed and for any continuous mapping $f : X\longrightarrow Y$
 to a  $\mathbf{K}$-space $Y$, there exists a unique $y_A\in Y$ such that $\overline{f(A)}=\overline{\{y_A\}}\}$. Endow $\mathbf{K}(X)$ with the lower Vietoris topology and denote the resulting space by $P_H(\mathbf{K}(X))$. We prove that if $P_H(\mathbf{K}(X))$ is a $\mathbf{K}$-space,  then the pair $\langle X^k=P_H(\mathbf{K}(X)), \eta_X\rangle$, where $\eta_X :X\longrightarrow X^k$, $x\mapsto\overline{\{x\}}$, is the $\mathbf{K}$-reflection of $X$. We call $\mathbf{K}$ an adequate category if for any $T_0$ space $X$, $P_H(\mathbf{K}(X))$ is a $\mathbf{K}$-space. So if $\mathbf{K}$ is adequate, then $\mathbf{K}$ is reflective in $\mathbf{Top}_0$. We show that $\mathbf{Sob}$, $\mathbf{Top}_d$, $\mathbf{Top}_w$ and the Keimel and Lawson's category are all adequate. Therefore, they are all reflective in $\mathbf{Top}_0$. Some major properties of $\mathbf{K}$-spaces and $\mathbf{K}$-reflections of $T_0$ spaces are investigated. In particular, it is proved that if $\mathbf{K}$ is adequate, then the $\mathbf{K}$-reflection preserves finite products of $T_0$ spaces. More precisely, for a finitely family $\{X_i: 1\leq i\leq n\}$ of $T_0$ spaces, we have that $(\prod\limits_{i=1}^{n}X_i)^k=\prod\limits_{i=1}^{n}X_i^k$ (up to homeomorphism). Our study also leads to a number of problems, whose answering will deepen our understanding of the related spaces and their categorical structures.

\section{Preliminary}

In this section, we briefly  recall some fundamental concepts and notations that will be used in the paper. Some basic properties of irreducible sets are presented.

For a poset $P$ and $A\subseteq P$, let
$\mathord{\downarrow}A=\{x\in P: x\leq  a \mbox{ for some }
a\in A\}$ and $\mathord{\uparrow}A=\{x\in P: x\geq  a \mbox{
	for some } a\in A\}$. For  $x\in P$, we write
$\mathord{\downarrow}x$ for $\mathord{\downarrow}\{x\}$ and
$\mathord{\uparrow}x$ for $\mathord{\uparrow}\{x\}$. A subset $A$
is called a \emph{lower set} (resp., an \emph{upper set}) if
$A=\mathord{\downarrow}A$ (resp., $A=\mathord{\uparrow}A$). Let $P^{(<\omega)}=\{F\subseteq P : F \mbox{~is a nonempty finite set}\}$ and $\mathbf{Fin} ~P=\{\uparrow F : F\in P^{(<\omega)}\}$.

The category of all $T_0$ spaces is denoted by $\mathbf{Top}_0$. For $X\in \mathbf{Top}_0$, we use $\leq_X$ to represent the \emph{specialization order} of $X$, that is, $x\leq_X y$ if{}f $x\in \overline{\{y\}}$). In this paper, when a $T_0$ space $X$ is considered as a poset, the order always refers to the specialization order if no other explanation. Let $\mathcal O(X)$ (resp., $\mathcal C(X)$) be the set of all open subsets (resp., closed subsets) of $X$, and let $\mathcal S^u(X)=\{\ua x : x\in X\}$. Define $\mathcal S_c(X)=\{\overline{{\{x\}}} : x\in X\}$ and $\mathcal D_c(X)=\{\overline{D} : D\in \mathcal D(X)\}$.

A nonempty subset $D$ of a poset $P$ is \emph{directed} if every two
elements in $D$ have an upper bound in $D$. The set of all directed sets of $P$ is denoted by $\mathcal D(P)$.   $P$ is called a
\emph{directed complete poset}, or \emph{dcpo} for short, if for any
$D\in \mathcal D(P)$, $\bigvee D$ exists in $P$. A subset $U$ of $P$ is \emph{Scott open} if
(i) $U=\mathord{\uparrow}U$ and (ii) for any directed subset $D$ for
which $\bigvee D$ exists, $\bigvee D\in U$ implies $D\cap
U\neq\emptyset$. All Scott open subsets of $P$ form a topology,
and we call this topology  the \emph{Scott topology} on $P$ and
denote it by $\sigma(P)$. The space $\Sigma~\!\! P=(P,\sigma(P))$ is called the
\emph{Scott space} of $P$.

A $T_0$ space $X$ is called a $d$-\emph{space} (or \emph{monotone convergence space}) if $X$ (with the specialization order) is a dcpo
 and $\mathcal O(X) \subseteq \sigma(X)$ (cf. \cite{redbook, Wyler}). Clearly, for a dcpo $P$, $\Sigma~\!\! P$ is a $d$-space. The category of all $d$-spaces with continuous mappings is denoted by $\mathbf{Top}_d$.

One can directly get the following result (cf. \cite{xu-shen-xi-zhao-19}).

\begin{proposition}\label{d-spacecharac1} For a $T_0$ space $X$, the following conditions are equivalent:
\begin{enumerate}[\rm (1) ]
	        \item $X$ is a $d$-space.
            \item $\mathcal D_c(X)=\mathcal S_c(X)$.
\end{enumerate}
\end{proposition}

The following result is well-known (cf. \cite{redbook}).

\begin{lemma}\label{Scott-cont1} Let $P, Q$ be posets and $f : P \longrightarrow Q$. Then the following two conditions are equivalent:
\begin{enumerate}[\rm (1)]
	\item $f$ is Scott continuous, that is, $f : \Sigma~\!\! P \longrightarrow \Sigma~\!\! Q$ is continuous.
	\item For any $D\in \mathcal D(P)$ for which $\bigvee D$ exists, $f(\bigvee D)=\bigvee f(D)$.
\end{enumerate}
\end{lemma}

\begin{lemma}\label{KL-basiclemma}\emph{(\cite{Keimel-Lawson})} Let $f : X \longrightarrow Y$ be a continuous mapping of $T_0$ spaces. If $D\in \mathcal {D}(X)$ has a supremum to which it converges, then $f(\bigvee D)=\bigvee f(D)$.
\end{lemma}

\begin{corollary}\label{continuous-ScottCONT} Let $P$ be a poset and $Y$ a $T_0$ space. If $f : \Sigma~\!\! P \longrightarrow Y$ is continuous, then $f : \Sigma~\!\! P \longrightarrow \Sigma~\!\! Y$ is continuous.
\end{corollary}

Let $\mathbf{Poset}$ denote the category of all posets with monotone (i.e. order-preserving) mappings, $\mathbf{DCPO}$ the category of all dcpos with Scott continuous mappings, and $\mathbf{Poset}_s$ the category of all posets with Scott continuous mappings. Then $\mathbf{DCPO}$ is a full subcategory of $\mathbf{Poset}_s$.

 A nonempty subset $A$ of  a $T_0$ space $X$ is \emph{irreducible} if for any $\{F_1, F_2\}\subseteq \mathcal C(X)$, $A \subseteq F_1\cup F_2$ implies $A \subseteq F_1$ or $A \subseteq  F_2$.  Denote by $\ir(X)$ (resp., $\ir_c(X)$) the set of all irreducible (resp., irreducible closed) subsets of $X$. Clearly, every subset of $X$ that is directed under $\leq_X$ is irreducible. $X$ is called \emph{sober}, if for any  $F\in\ir_c(X)$, there is a unique point $a\in X$ such that $F=\overline{\{a\}}$. The category of all sober spaces with continuous mappings is denoted by $\mathbf{Sob}$.

The following two lemmas on irreducible sets are well-known.

\begin{lemma}\label{irrsubspace}
Let $X$ be a space and $Y$ a subspace of $X$. Then the following conditions are equivalent for a
subset $A\subseteq Y$:
\begin{enumerate}[\rm (1)]
	\item $A$ is an irreducible subset of $Y$.
	\item $A$ is an irreducible subset of $X$.
	\item ${\rm cl}_X A$ is an irreducible closed subset of $X$.
\end{enumerate}
\end{lemma}

\begin{lemma}\label{irrimage}
	If $f : X \longrightarrow Y$ is continuous and $A\in\ir (X)$, then $f(A)\in \ir (Y)$.
\end{lemma}

\begin{lemma}\label{irrprod}\emph{(\cite{Shenchon})}
	Let	$X=\prod_{i\in I}X_i$ be the product space of $T_0$ spaces  $X_i (i\in I)$. If  $A$ is an irreducible subset of $X$, then $\cl_X(A)=\prod_{i\in I}\cl_{X_i}(p_i(A))$, where $p_i : X \longrightarrow X_i$ is the $i$th projection for each $i\in I$.
\end{lemma}

\begin{lemma}\label{prodirr}\emph{(\cite{xu-shen-xi-zhao-19})}
	Let	$X=\prod_{i\in I}X_i$ be the product space of $T_0$ spaces  $X_i (i\in I)$  and $A_i\subseteq X_i$ for each $i\in I$. Then the following two conditions are equivalent:
\begin{enumerate}[\rm (1)]
	\item $\prod_{i\in I}A_i\in \ir (X)$.
	\item $A_i\in \ir (X_i)$ for each $i\in I$.
\end{enumerate}
\end{lemma}

By Lemma \ref{irrprod} and Lemma \ref{prodirr}, we obtain the following corollary.

\begin{corollary}\label{irrcprod} Let $X=\prod_{i\in I}X_i$ be the product space of $T_0$ spaces  $X_i (i\in I)$. If  $A\in \ir_c(X)$, then $A=\prod_{i\in I}p_i(A)$ and $p_i(A)\in \ir_c (X_i)$ for each $i\in I$.
\end{corollary}

For any topological space $X$, $\mathcal G\subseteq 2^{X}$ and $A\subseteq X$, let $\Diamond_{\mathcal G} A=\{G\in \mathcal G : G\bigcap A\neq\emptyset\}$ and $\Box_{\mathcal G} A=\{G\in \mathcal G : G\subseteq  A\}$. The symbols $\Diamond_{\mathcal G} A$ and $\Box_{\mathcal G} A$ will be simply written as $\Diamond A$  and $\Box A$ respectively if there is no confusion. The \emph{lower Vietoris topology} on $\mathcal{G}$ is the topology that has $\{\Diamond U : U\in \mathcal O(X)\}$ as a subbase, and the resulting space is denoted by $P_H(\mathcal{G})$. If $\mathcal{G}\subseteq \ir (X)$, then $\{\Diamond_{\mathcal{G}} U : U\in \mathcal O(X)\}$ is a topology on $\mathcal{G}$.

\begin{remark} \label{eta continuous} Let $X$ be a $T_0$ space.
\begin{enumerate}[\rm (1)]
	\item If $\mathcal{S}_c(X)\subseteq \mathcal{G}$, then the specialization order on $P_H(\mathcal{G})$ is the order of set inclusion, and the \emph{canonical mapping} $\eta_{X}: X\longrightarrow P_H(\mathcal{G})$, given by $\eta_X(x)=\overline {\{x\}}$, is an order and topological embedding (cf. \cite{redbook, Jean-2013, Schalk}).
    \item The space $X^s=P_H(\ir_c(X))$ with the canonical mapping $\eta_{X}: X\longrightarrow X^s$ is the \emph{sobrification} of $X$
    (cf. \cite{redbook, Jean-2013}).
        \item $P_H(\mathcal S_c(X))$ is a subspace of $X^s$ and $X$ is homeomorphic to $P_S(\mathcal S_c(X))$ via a homeomorphism $x\mapsto \overline{\{x\}}$.
\end{enumerate}
\end{remark}

\begin{remark}\label{PH-subtopology} Let $X$ be a $T_0$ space and $\mathcal G_1\subseteq \mathcal G_2\subseteq\ir_c(X)$. If $\mathcal G_2$ is endowed with  the lower Vietoris topology, then the subspace  $(\mathcal G_1, \{\Diamond U \cap \mathcal G_1 : U\in \mathcal O(X)\})$ is the space  $P_H(\mathcal{G}_1)$, which is a subspace of $X^s$. In what follows, when a subset $\mathcal G$ in $X^s$ (that is $\mathcal G\subseteq \ir_c(X)$) is considered as a topological space, the topology always refers to the subspace topology of $X^s$ if no other explanation.
\end{remark}

For a space $X$, a subset $A$ of $X$ is called \emph{saturated} if $A$ equals the intersection of all open sets containing it (equivalently, $A$ is an upper set in the specialization order). We shall use $\mathord{\mathsf{Q}}(X)$ to
denote the set of all nonempty compact saturated subsets of $X$ and endow it with the \emph{Smyth preorder}, that is, for $K_1,K_2\in \mathord{\mathsf{Q}}(X)$, $K_1\sqsubseteq K_2$ if{}f $K_2\subseteq K_1$. $X$ is called \emph{well-filtered} if it is $T_0$, and for any open set $U$ and filtered family $\mathcal{K}\subseteq \mathord{\mathsf{Q}}(X)$, $\bigcap\mathcal{K}{\subseteq} U$ implies $K{\subseteq} U$ for some $K{\in}\mathcal{K}$. The category of all well-filtered spaces with continuous mappings is denoted by $\mathbf{Top}_w$.
The space $P_S(\mathord{\mathsf{Q}}(X))$, denoted shortly by $P_S(X)$, is called the \emph{Smyth power space} or \emph{upper space} of $X$ (cf. \cite{Heckmann, Schalk}). It is easy to see that the specialization order on $P_S(X)$ is the Smyth order (that is, $\leq_{P_S(X)}=\sqsubseteq$). The \emph{canonical mapping} $\xi_X: X\longrightarrow P_S(X)$, $x\mapsto\ua x$, is an order and topological embedding (cf. \cite{Heckmann, Klause-Heckmann, Schalk}). Clearly, $P_S(\mathcal S^u(X))$ is a subspace of $P_S(X)$ and $X$ is homeomorphic to $P_S(\mathcal S^u(X))$ via a homeomorphism $x\mapsto \ua x$.

As in \cite{E_20182}, a topological space $X$ is \emph{locally hypercompact} if for each $x\in X$ and each open neighborhood $U$ of $x$, there is  $\ua F\in \mathbf{Fin}~X$ such that $x\in\ii\,\ua F\subseteq\ua F\subseteq U$. A space $X$ is called a $C$-\emph{space} if for each $x\in X$ and each open neighborhood $U$ of $x$, there is $u\in X$ such that $x\in\ii\,\ua u\subseteq\ua u\subseteq U$. A set $K\subseteq X$ is called \emph{supercompact} if for
any arbitrary family $\{U_i : i\in I\}\subseteq \mathcal O(X)$, $K\subseteq \bigcup_{i\in I} U_i$  implies $K\subseteq U$ for some $i\in I$. It is easy to check that the supercompact saturated sets of $X$ are exactly the sets $\ua x$ with $x \in X$ (see \cite[Fact 2.2]{Klause-Heckmann}). It is well-known that $X$ is a $C$-space if{}f $\mathcal O(X)$ is a \emph{completely distributive} lattice (cf. \cite{E_2009}). A space $X$ is called \emph{core compact} if $\mathcal O(X)$ is a \emph{continuous lattice} (cf. \cite{redbook}).

\begin{theorem}\label{WFLC=CoreC}\emph{(\cite{Lawson-Xi, xu-shen-xi-zhao-19})} Let $X$ be a well-filtered space. Then $X$ is locally compact if{}f $X$ is core compact.
\end{theorem}

\section{$\mathbf{K}$-sets}

For a full subcategory $\mathbf{K}$ of $\mathbf{Top}_0$, the objects of $\mathbf{K}$ will be called $\mathbf{K}$-spaces. In \cite{Keimel-Lawson}, Keimel and Lawson required the following properties:

($\mathrm{K}_1$) Homeomorphic copies of $\mathbf{K}$-spaces are $\mathbf{K}$-spaces.

($\mathrm{K}_2$) All sober spaces are $\mathbf{K}$-spaces or, equivalently, $\mathbf{Sob}\subseteq \mathbf{K}$.

($\mathrm{K}_3$) In a sober space S, the intersection of any family of $\mathbf{K}$-subspaces is a $\mathbf{K}$-space.

($\mathrm{K}_4$) Continuous maps $f : S \longrightarrow T$ between sober spaces $S$ and $T$ are $\mathbf{K}$-continuous, that is, for every $\mathbf{K}$-subspace $K$ of $T$ , the inverse image $f^{-1}(K)$ is a $\mathbf{K}$-subspace of $S$.

\begin{definition}\label{KL-category} A full subcategory $\mathbf{K}$ of $\mathbf{Top}_0$ is said to be \emph{closed with respect to homeomorphisms} if $\mathbf{K}$ has ($\mathrm{K}_1$). $\mathbf{K}$ is called a \emph{Keimel-Lawson category} if $\mathbf{K}$ satisfies ($\mathrm{K}_1$)-($\mathrm{K}_4$).

\end{definition}

Clearly, $\mathbf{Sob}, \mathbf{Top}_d$ and $\mathbf{Top}_w$ are closed with respect to homeomorphisms and satisfy ($\mathrm{K}_2$).

In what follows, $\mathbf{K}$ always refers to a full subcategory $\mathbf{Top}_0$ containing $\mathbf{Sob}$, that is,  $\mathbf{K}$ has ($\mathrm{K}_2$). For two spaces $X$ and $Y$, we use the symbol $X\cong Y$ to  represent that $X$ and $Y$ are homeomorphic.

\begin{definition}\label{K subset}
	 A subset $A$ of a $T_0$ space $X$ is called a $\mathbf{K}$-\emph{set}, provided for any continuous mapping $ f:X\longrightarrow Y$
to a $\mathbf{K}$-space $Y$, there exists a unique $y_A\in Y$ such that $\overline{f(A)}=\overline{\{y_A\}}$.
Denote by $\mathbf{K}(X)$ the set of all closed $\mathbf{K}$-sets of $X$. $X$ is said to be a $\mathbf{K}$-\emph{determined} space if $\ir_c(X)=\mathbf{K}(X)$ or, equivalently, all irreducible closed sets of $X$ are $\mathbf{K}$-sets (it is easy to check that all $\mathbf{K}$-sets are irreducible, please see Corollary \ref{SKIsetrelation} below).
\end{definition}

Obviously, a subset $A$ of a space $X$ is a $\mathbf{K}$-set if{}f $\overline{A}$ is a $\mathbf{K}$-set. For simplicity, let $d(X)=\mathbf{Top}_d(X)$ and $\mathsf{WF}(X)=\mathsf{Top}_w(X)$. $X$ is called \emph{well-filtered determined}, $\mathsf{WF}$-determined for short, if all irreducible closed subsets of $X$ are $\mathsf{WF}$-sets, that is, $\ir_c(X)=\mathsf{WF} (X)$.

\begin{lemma}\label{sobd=irr} For a $T_0$ space $X$,  $\mathbf{Sob}(X)=\ir_c(X)$.
\end{lemma}

\begin{proof} By Lemma \ref{irrimage}, $\ir_c(X)\subseteq \mathbf{Sob}(X)$. Suppose $A\in \mathbf{Sob}(X)$. Now we show that $A$ is irreducible. Consider the sobrification $X^s$ ($=P_H(\ir_c(X)$) of $X$ and the canonical topological embedding $\eta_{X}: X\longrightarrow X^s$, given by $\eta_X(x)=\overline {\{x\}}$. Then there is a $B\in \ir_c(X)$ such that $\Box_{\ir_c(X)}A=\overline{\eta_X (A)}=\overline {\{B\}}=\Box_{\ir_c(X)}B$, and whence $A=B$. Thus $A\in \ir_c(X)$.
\end{proof}

\begin{corollary}\label{SKIsetrelation} For a $T_0$ space $X$,  $\mathcal S_c(X)\subseteq\mathbf{K}(X)\subseteq\ir_c(X)$.
\end{corollary}
\begin{proof}
	Clearly, $\mathcal{S}_c(X)\subseteq \mathbf{K}(X)$. Since $\mathbf{Sob}\subseteq \mathbf{K}$, we have $\mathbf{K}(X)\subseteq\mathbf{Sob}(X)=\ir_c(X)$ by Lemma \ref{sobd=irr}.
\end{proof}

Rudin's Lemma \cite{Rudin} is a very useful tool in domain theory and non-Hausdorff topology (see [3, 6-9, 11, 20, 25, 26]). In \cite{Klause-Heckmann}, Heckman and Keimel presented the following topological variant of Rudin's Lemma.

\begin{lemma}\label{t Rudin} \emph{(Topological Rudin's Lemma)} Let $X$ be a topological space and $\mathcal{A}$ an
irreducible subset of the Smyth power space $P_S(X)$. Then every closed set $C {\subseteq} X$  that
meets all members of $\mathcal{A}$ contains an minimal irreducible closed subset $A$ that still meets all
members of $\mathcal{A}$.
\end{lemma}

For a $T_0$ space $X$ and $\mathcal{K}\subseteq \mathord{\mathsf{Q}}(X)$, let $M(\mathcal{K})=\{A\in \mathcal C(X) : K\bigcap A\neq\emptyset \mbox{~for all~} K\in \mathcal{K}\}$ (that is, $\mathcal A\subseteq \Diamond A$) and $m(\mathcal{K})=\{A\in \mathcal C(X) : A \mbox{~is a minimal menber of~} M(\mathcal{K})\}$. The following concept was introduced based on topological Rudin's Lemma.

\begin{definition}\label{rudinset} (\cite{Shenchon, xu-shen-xi-zhao-19})
		Let $X$ be a $T_0$ space. A nonempty subset  $A$  of $X$  is said to have the \emph{Rudin property}, if there exists a filtered family $\mathcal K\subseteq \mathord{\mathsf{Q}}(X)$ such that $\overline{A}\in m(\mathcal K)$ (that is,  $\overline{A}$ is a minimal closed set that intersects all members of $\mathcal K$). Let $\mathsf{RD}(X)=\{A\in \mathcal C(X) : A\mbox{~has Rudin property}\}$. The sets in $\mathsf{RD}(X)$ will also be called \emph{Rudin sets}.
\end{definition}

\begin{lemma}\label{rudinwf}
	Let $X$ be a $T_0$ space and  $Y$ a well-filtered space. If $f : X\longrightarrow Y$ is continuous and $A\subseteq X$ has Rudin property, then there exists a unique $y_A\in X$ such that $\overline{f(A)}=\overline{\{y_A\}}$.
\end{lemma}
\begin{proof}
Since $A$ has Rudin property, there exists a filtered family $\mathcal K\subseteq \mathord{\mathsf{Q}}(X)$ such that $\overline{A}\in m(\mathcal K)$. Let $\mathcal{K}_f=\{\ua f(K\cap \overline{A}) : K\in \mathcal K\}$. Then $\mathcal{F}_f\subseteq \mathord{\mathsf{Q}}(Y)$ is filtered. For each $K\in \mathcal{K}$, since $K\cap A\neq\emptyset$, we have $\emptyset\neq f(K\cap A)\subseteq \ua f(K\cap A)\cap \overline{f(A)}$. So $\overline{f(A)}\in M(\mathcal{K}_f)$. If $B$ is a closed subset of $\overline{f(A)}$ with $B\in M(\mathcal{K}_f)$, then $B\cap\ua f(K\cap A)\neq\emptyset$ for every $K\in \mathcal K$. So $K\cap A\cap f^{-1}(B)\neq\emptyset$ for all $K\in \mathcal K$. It follows that $A=A\cap f^{-1}(B)$ by the minimality of $A$, and consequently, $\overline{f(A)}\subseteq B$. Therefore, $\overline{f(A)}=B$. Thus $\overline{f(A)}\in m(\mathcal{K}_f)$. Since $Y$ is well-filtered, we have $\bigcap_{K\in \mathcal{K}}\ua f(K\cap \overline{A})\cap \overline{f(A)}\neq\emptyset$. Select a $y_A\in \bigcap_{K\in \mathcal K} \ua f(K\cap \overline{A})\cap \overline{f(A)}$. Then $\overline{\{y_A\}}\subseteq \overline{f(A)}$ and $K\cap \overline{A}\cap f^{-1}(\overline{\{y_A\}})\neq\emptyset$ for all $K\in \mathcal K$. It follows that $\overline{A}=\overline{A}\cap f^{-1}(\overline{\{y_A\}})$ by the minimality of $\overline{A}$, and consequently, $\overline{f(A)}\subseteq \overline{\{y_A\}}$. Therefore, $\overline{f(A)}=\overline{\{y_A\}}$. The uniqueness of $y_A$ follows from the $T_0$ separation of $Y$.
\end{proof}

\begin{proposition}\label{DRWIsetrelation}
	Let $X$ be a $T_0$ space. Then $\mathcal{S}_c(X)\subseteq\mathcal{D}_c(X)\subseteq \mathsf{RD}(X)\subseteq\mathsf{WF}(X)\subseteq\ir_c(X)$.
\end{proposition}
\begin{proof} Obviously, $\mathcal{S}_c(X)\subseteq \mathcal{D}_c$. Now we prove that the closure of a directed subset $D$ of $X$ is a Rudin set. Let
$\mathcal K_D=\{\ua d : d\in D\}$. Then $\mathcal K_D\subseteq \mathord{\mathsf{Q}}(X)$ is filtered and $\overline{D}\in M(\mathcal K_D)$. If $A\in M(\mathcal K_D)$, then $d\in A$ for every $d\in D$, and hence $\overline{D}\subseteq A$. So $\overline{D}\in m(\mathcal K_D)$. Therefore $\overline{D}\in \mathsf{RD}(X)$. By Lemma \ref{rudinwf}, $\mathsf{RD}(X)\subseteq\mathsf{WF}(X)$. Finally, by Corollary \ref{SKIsetrelation} (for $\mathbf{K}=\mathbf{Top}_w$), we have $\mathsf{WF}(X)\subseteq\ir_c(X)$.
\end{proof}

\begin{example}\label{K-examp1}
	Let $X$ be a countable infinite set and endow $X$ with the cofinite topology (having the complements of the finite sets as open sets). The
resulting space is denoted by $X_{cof}$. Then $\mathsf{Q} (X_{cof})=2^X\setminus \{\emptyset\}$ (that is, all nonempty subsets of $X$), and hence $X_{cof}$ is a locally compact and first countable $T_1$ space.  Let $\mathcal K=\{X\setminus F : F\in X^{(<\omega)}\}$. It is easy to check that $\mathcal K\subseteq \mathsf{Q} (X_{cof})$ is filtered and $X\in m(\mathcal K)$. Therefore, $X\in \kf(X)$ but $X\not\in \md_c(X)$, and whence $\kf(X)\neq \md_c(X)$ and $\mathsf{WF} (X)\neq \md_c(X)$. Thus $X_{cof}$ is not well-filtered (and hence non-sober).
\end{example}

\begin{example}\label{K-examp2}
	Let $L$ be the complete lattice constructed by Isbell \cite{isbell} and $\mathbf{K}=\mathbf{Top}_w$. Then by \cite[Corollary 3.2]{Xi-Lawson-2017}, $\Sigma L$ is a well-filtered space, and whence $\mathsf{WF}(X)=\mathcal S_c(X)$. Note that $\Sigma L$ is not sober. Therefore, by Prpposition \ref{DRWIsetrelation}, $\mathsf{WF}(X)\neq\ir_c(X)$ and $\kf(X)\neq\ir_c(X)$.
\end{example}

\begin{lemma}\label{K-setimage}
Let $X,Y$ be two $T_0$ spaces. If $f:X\longrightarrow Y$ is a continuous mapping and $A\in \mathbf{K} (X)$, then $\overline{f(A)}\in \mathbf{K} (Y)$.
\end{lemma}
\begin{proof}	Let $Z$ is a $\mathbf{K}$-space and $g:Y\longrightarrow Z$ is a continuous mapping.
Since $g\circ f:X\longrightarrow Z$ is continuous and $A\in \mathbf{K} (X)$, there is $z\in Z$ such that $\overline{g(\overline{f(A)})}=\overline{g\circ f(A)}=\overline{\{z\}}$. Thus $\overline{f(A)}\in \mathbf{K} (Y)$.
\end{proof}

\begin{lemma}\label{K-setprod}
	Let	$\{X_i: 1\leq i\leq n\}$ be a finite family of $T_0$ spaces and $X=\prod\limits_{i=1}^{n}X_i$ the product space. For $A\in\ir (X)$, the following conditions are equivalent:
\begin{enumerate}[\rm (1)]
	\item $A$ is a $\mathbf{K}$-set.
	\item $p_i(A)$ is a $\mathbf{K}$-set for each $1\leq i\leq n$.
\end{enumerate}
\end{lemma}

\begin{proof} (1) $\Rightarrow$ (2): By Lemma \ref{K-setimage}.

(2) $\Rightarrow$ (1): By induction, we need only to prove the implication for the case of $n=2$. Let $A_1=\cl_{X_1} p_1(A)$ and $A_2=\cl_{X_2} p_2(A)$. Then by condition (2), $(A_1, A_2)\in \mathbf{K} (X_1)\times \mathbf{K} (X_2)$. Now we show that the product $A_1\times A_2\in\mathbf{K} (X)$. Let $f : X_1\times X_2 \longrightarrow Y$ a continuous mapping from $X_1\times X_2$ to a $\mathbf{K}$-space $Y$. For each $b\in X_2$, $X_1$ is homeomorphic to $X_1\times \{b\}$ (as a subspace of $X_1\times X_2$) via the homeomorphism $\mu_b : X_1 \longrightarrow X_1\times \{b\}$ defined by $\mu_b(x)=(x, b)$. Let $i_b : X_1\times \{b\}\longrightarrow X_1\times X_2$ be the embedding of $X_1\times \{b\}$ in $X_1\times X_2$. Then  $f_{b}=f\circ i_b \circ \mu_b : X_1 \longrightarrow Y$, $f_b(x)=f((x, b))$,  is continuous. Since $A_1\in \mathbf{K} (X_1)$, there is a unique $y_b\in Y$ such that $\overline{f(A_1\times \{b\})}=\overline{f_b(A_1)}=\overline{\{y_b\}}$. Define a mapping $g_A : X_2 \longrightarrow Y$ by $g_A(b)=y_b$. For each $V\in \mathcal O(Y)$,
$$\begin{array}{lll}
	g_A^{-1}(V)& =\{b\in X_2 : g_A(b)\in V\}\\
	           & =\{b\in X_2 : \overline{f_b(A_1)}\cap V\neq\emptyset\}\\
	           & =\{b\in X_2 : \overline{f(A_1\times \{b\})}\cap V\neq\emptyset\}\\
	           & =\{b\in X_2 : f(A_1\times \{b\})\cap V\neq\emptyset\}\\
               & =\{b\in X_2 : (A_1\times \{b\})\cap f^{-1}(V)\neq\emptyset\}.\\
	\end{array}$$
 Therefore, for each $b\in g_A^{-1}(V)$, there is an $a_1\in A_1$ such that $(a_1, b)\in f^{-1}(V)\in \mathcal O(X_1\times X_2)$, and hence there is $(U_1, U_2)\in \mathcal O(X_1)\times \mathcal O(X_2)$ such that $(a_1, b)\in U_1\times U_2\subseteq  f^{-1}(V)$. It follows that $b\in U_2\subseteq g_A^{-1}(V)$. Thus $g_A : X_2 \longrightarrow Y$ is continuous. Since $A_2\in \mathbf{K} (X_2)$, there is a unique $y_A\in Y$ such that $\overline{g_A(A_2)}=\overline{\{y_A\}}$. Therefore, by Lemma \ref{irrprod}, we have
 $$\begin{array}{lll}
     \overline{f(\cl_X A)} & =\overline{f(A_1\times A_2)}\\
	                       & =\overline{\bigcup\limits_{a_2\in A_2}f(A_1\times \{a_2\})}\\
	                       & =\overline{\bigcup\limits_{a_2\in A_2}\overline{f(A_1\times \{a_2\})}}\\
	                       & =\overline{\bigcup\limits_{a_2\in A_2}\overline{\{g_A(a_2)\}}}\\
                           & =\overline{\bigcup\limits_{a_2\in A_2}\{g_A(a_2)\}}\\
                           & =\overline{g_A(A_2)}\\
                           & =\overline{\{y_A\}}.\\
	\end{array}$$
Thus $\cl_X A\in \mathbf{K}(X)$, and hence $A$ is a $\mathbf{K}$-set.
\end{proof}

By Corollary \ref{irrcprod} and Lemma \ref{K-setprod}, we get the following result.

\begin{corollary}\label{K-closedsetprod}
	Let	$X=\prod\limits_{i=1}^{n}X_i$ be the product of a finitely family $\{X_i: 1\leq i\leq n\}$ of $T_0$ spaces. If $A\in\mathbf{K} (X)$, then $A=\prod\limits_{i=1}^{n}p_i(X_i)$, and $p_i(A)\in \mathbf{K} (X_i)$ for all $1\leq i \leq n$.
\end{corollary}

\section{A direct construction of $K$-reflections of $T_0$ spaces}

In Section 4, we give a direct construction of the  $\mathbf{K}$-reflections of $T_0$ spaces and investigate some basic properties of $\mathbf{K}$-spaces and  $\mathbf{K}$-reflections. In particular, it is proved that if $\mathbf{K}$ ia an adequate category, then the $\mathbf{K}$-reflection preserves finite products of $T_0$ spaces.

\begin{definition}\label{WFtion}
	Let $X$ be a $T_0$ space. A $\mathbf{K}$-\emph{reflection} of $X$ is a pair $\langle \widetilde{X}, \mu\rangle$ consisting of a $\mathbf{K}$-space $\widetilde{X}$ and a continuous mapping $\mu :X\longrightarrow \widetilde{X}$ satisfying that for any continuous mapping $f: X\longrightarrow Y$ to a $\mathbf{K}$-space, there exists a unique continuous mapping $f^* : \widetilde{X}\longrightarrow Y$ such that $f^*\circ\mu=f$, that is, the following diagram commutes.\\
\begin{equation*}
	\xymatrix{
		X \ar[dr]_-{f} \ar[r]^-{\mu}
		&\widetilde{X}\ar@{.>}[d]^-{f^*}\\
		&Y}
	\end{equation*}

\end{definition}

By a standard argument, $\mathbf{K}$-reflections, if they exist, are unique up to homeomorphism. We shall use $X^k$ to denote the space of the $\mathbf{K}$-reflection of $X$ if it exists.

By Corollary \ref{SKIsetrelation}, $\{\Diamond_{\mathbf{K}(X)} U : U\in \mathcal O(X)\}$ is a topology on $\mathbf{K}(X)$. In the following, let $\eta_X^k : X\longrightarrow P_H(\mathbf{K}(X))$, $\eta_X^k(x)=\overline {\{x\}}$, be the canonical mapping from $X$ to $P_H(\mathbf{K}(X))$.

\begin{lemma}\label{lemmaeta}
	The canonical mapping $\eta_X^k:X\longrightarrow P_H(\mathbf{K}(X))$ is a topological embedding.
\end{lemma}
\begin{proof}
For $U\in\mathcal O(X)$, we have $$(\eta_X^{k})^{-1}(\Diamond U)=\{x\in X: \da x\in\Diamond U\}=\{x\in X: x\in U\}=U,$$ so $\eta_X^k$ is continuous.
In addition, we have $$
\eta_X^k(U)=\{\da x: x\in U\}
=\{\da x: \da x\in\Diamond U\}
=\Diamond U\cap \eta_X^k(X),$$
which implies that $\eta_X^k$ is an open mapping to $\eta_X^k(X)$, as a subspace of $P_H(\mathbf{K}(X))$.
As $\eta_X^k$ is injective, $\eta_X^k$ is a topological embedding.
\end{proof}

\begin{lemma}\label{lemmaclosure}
	For a $T_0$ space $X$ be and $A\subseteq X$, $\overline{\eta_X^k(A)}=\overline{\eta_X^k\left(\overline{A}\right)}=\overline{\Box A}=\Box \overline{A}$ in $P_H(\mathbf{K}(X))$.
\end{lemma}
\begin{proof}
	Clearly, $\eta_X^k(A)\subseteq \Box A\subseteq \Box\overline{A}$, $\eta_X^k\left(\overline{A}\right)\subseteq \Box\overline{A}$ and $\Box\overline{A}$ is closed in $P_H(\mathbf{K}(X))$. It follows that
	$$\overline{\eta_X^k(A)}\subseteq \overline{\Box A}\subseteq \Box \overline{A}\ \text{ and }\ \overline{\eta_X^k(A)}\subseteq\overline{\eta_X^k\left(\overline{A}\right)}\subseteq\Box\overline{A}.$$
	To complete the proof, we need to show $\Box\overline{A}\subseteq \overline{\eta_X^k(A)}$.
	Let $F\in \Box\overline{A}$. Suppose $U\in\mathcal O(X)$ such that $F\in\Diamond U$ (note that $F\in \mathbf{K}(X)$), that is, $F\cap U\neq\emptyset$. Since $F\subseteq \overline{A}$, we have $A\cap U\neq\emptyset$. Let $a\in A\cap U$. Then $\da a\in \Diamond U\cap \eta_X^k(A)\neq\emptyset$. This implies that $F\in \overline{\eta_X^k(A)}$. Thus $\Box\overline{A}\subseteq \overline{\eta_X^k(A)}$.	
\end{proof}

\begin{lemma}\label{lemmaK-irr}
	Let $X$ be a $T_0$ space and $A$ a nonempty subset of $X$. Then the following conditions are equivalent:
	\begin{enumerate}[\rm (1)]
		\item $A$ is irreducible in $X$.
		\item $\Box A$ is irreducible in $P_H(\mathbf{K}(X))$.
        \item $\Box \overline{A}$ is irreducible in $P_H(\mathbf{K}(X))$.
	\end{enumerate}
\end{lemma}
\begin{proof}
	(1) $\Rightarrow$ (3): Assume $A$ is irreducible. Then $\eta_X^k(A)$ is irreducible in $P_H(\mathbf{K}(X))$ by Lemma \ref{irrimage} and Lemma \ref{lemmaeta}. By Lemma \ref{irrsubspace} and Lemma \ref{lemmaclosure}, $\Box \overline{A}=\overline{\eta_X^k(A)}$ is irreducible in $P_H(\mathbf{K}(X))$.
	
	(3) $\Rightarrow$ (1): Assume $\Box \overline{A}$ is irreducible. Let $A\subseteq B\cup C$ with  $B,C\in\mathcal C(X)$. By Corollary \ref{SKIsetrelation}, $\mathbf{K} (X)\subseteq \ir_c(X)$, and consequently, we have $\Box\overline{A}\subseteq \Box B\cup\Box C$. Since $\Box\overline{A}$ is irreducible,  $\Box\overline{A}\subseteq \Box B$ or $\Box\overline{A}\subseteq C$, showing that  $\overline{A}\subseteq B$ or $\overline{A}\subseteq C$, and consequently, $A\subseteq B$ or $A\subseteq C$, proving $A$ is irreducible.

    (2) $\Leftrightarrow$ (3): By Lemma \ref{irrsubspace} and Lemma \ref{lemmaclosure}.

\end{proof}

For the $\mathbf{K}$-reflections of $T_0$ spaces, the following lemma is crucial.

\begin{lemma}\label{K-lemmafstar}
Let $X$ be a $T_0$ space and $f:X\longrightarrow Y$ a continuous mapping from $X$ to a well-filtered space $Y$. Then there exists a unique continuous mapping $f^* :P_H(\mathbf{K}(X))\longrightarrow Y$ such that $f^*\circ\eta_X^k=f$, that is, the following diagram commutes.
\begin{equation*}
\xymatrix{
	X \ar[dr]_-{f} \ar[r]^-{\eta_X^k}
	&P_H(\mathbf{K}(X))\ar@{.>}[d]^-{f^*}\\
	&Y}
\end{equation*}	
\end{lemma}
\begin{proof}For each $A\in\mathbf{K}(X)$, there exists a unique $y_A\in Y$ such that $\overline{f(A)}=\overline{\{y_A\}}$. Then we can define a mapping $f^*:P_H(\mathbf{K}(X))\longrightarrow Y$ by
$$\forall A\in\mathbf{K}(X),\ \ f^*(A)=y_A.$$

{Claim 1:}  $f^*\circ \eta_X^k=f$.

Let $x\in X$. Since $f$ is continuous, we have
$\overline{f\left(\overline{\{x\}}\right)}=\overline{f(\{x\})}=\overline{\{f(x)\}}$,
so
$f^*\left(\overline{\{x\}}\right)=f(x)$. Thus $f^*\circ \eta_X^k=f$.

{Claim 2:}  $f^*$ is continuous.

Let $V\in\mathcal O(Y)$. Then
$$\begin{array}{lll}
(f^*)^{-1}(V)&=&\{A\in\mathbf{K}(X): f^*(A)\in V\}\\
&=&\{A\in\mathbf{K}(X): \overline{\{f^*(A)\}}\cap V\neq\emptyset\}\\
&=&\{A\in\mathbf{K}(X): \overline{f(A)}\cap V\neq\emptyset\}\\
&=&\{A\in\mathbf{K}(X): f(A)\cap V\neq\emptyset\}\\
&=&\{A\in\mathbf{K}(X): A\cap f^{-1}(V)\neq\emptyset\}\\
&=&\Diamond f^{-1}(V),
\end{array}$$
which shows that $(f^*)^{-1}(V)$ is open in $P_H(\mathbf{K}(X))$. Thus  $f^*$ is continuous.

{Claim 3:}  The mapping $f^*$ is unique such that $f^*\circ \eta_X^k=f$.

Assume $g:P_H(\mathbf{K}(X))\longrightarrow Y$ is a continuous mapping such that $g\circ\eta_X^k=f$.  Let $A\in\mathbf{K}(X)$. We need to show $g(A)=f^*(A)$.
Let $a\in A$.  Then $\overline{\{a\}}\subseteq A$, implying that $g(\overline{\{a\}})\leq_Y g(A)$, that is,  $g\left(\overline{\{a\}}\right)=f(a)\in\overline{\{g(A)\}}$. Thus $\overline{\{f^*(A)\}}=\overline{f(A)}\subseteq \overline{\{g(A)\}}$.
In addition, since $A\in\overline{\eta_X^k(A)}$ and $g$ is continuous, $g(A)\in g\left(\overline{\eta_X^k(A)}\right)\subseteq\overline{g(\eta_X^k(A))}=\overline{f(A)}=\overline{\{f^*(A)\}}$, which implies that $\overline{\{g(A)\}}\subseteq \overline{\{f^*(A)\}}$.  So $\overline{\{g(A)\}}=\overline{\{f^*(A)\}}$. Since $Y$ is $T_0$, $g(A)=f^*(A)$. Thus $g=f^*$.
\end{proof}

From Lemma \ref{K-lemmafstar} we deduce the following main result of this paper.

\begin{theorem}\label{K-reflection}
	Let $X$ be a $T_0$ space. If $P_H(\mathbf{K}(X))$ is a $\mathbf{K}$-space,  then the pair $\langle X^k=P_H(\mathbf{K}(X)), \eta_X^k\rangle$, where $\eta_X^k :X\longrightarrow X^k$, $x\mapsto\overline{\{x\}}$, is the $\mathbf{K}$-reflection of $X$.
\end{theorem}

\begin{definition}\label{K-adequate} $\mathbf{K}$ is called \emph{adequate} if for any $T_0$ space $X$, $P_H(\mathbf{K}(X))$ is a $\mathbf{K}$-space.
\end{definition}

\begin{corollary}\label{K-adequate reflective}
	If $\mathbf{K}$ is adequate, then $\mathbf{K}$ is reflective in $\mathbf{Top}_0$.
\end{corollary}

\begin{corollary}\label{K-fuctor}
	If $\mathbf{K}$ is adequate, then for any  $T_0$ spaces $X,Y$ and any continuous mapping $f:X\longrightarrow Y$, there exists a unique continuous mapping $f^k:X^k\longrightarrow Y^k$ such that $f^k\circ \eta_X^k=\eta_Y^k\circ f$, that is, the following diagram commutes.
		\begin{equation*}
	\xymatrix{
		X \ar[d]_-{f} \ar[r]^-{\eta_X^k} &X^k\ar[d]^-{f^k}\\
		Y \ar[r]^-{\eta_Y^k} &Y^k
	}
	\end{equation*}
For each $A\in \mathbf{K} (X)$, $f^k(A)=\overline{f(A)}$.
\end{corollary}

Corollary \ref{K-fuctor} defines a functor $K : \mathbf{Top}_0 \longrightarrow \mathbf{K}$, which is the left adjoint to the inclusion functor $I : \mathbf{K} \longrightarrow \mathbf{Top}_0$.

\begin{corollary}\label{KKc}
	Suppose that $\mathbf{K}$ is adequate and closed with respect to homeomorphisms. Then for any $T_0$ space $X$, the following conditions are equivalent:
	\begin{enumerate}[\rm (1)]
		\item $X$ is a $\mathbf{K}$-space.
        \item $\mathbf{K} (X)=\mathcal S_c(X)$, that is, for each $A\in\mathbf{K}(X)$, there exists an $x\in X$ such that $A=\overline{\{x\}}$.
        \item $X\cong X^k$.
        \end{enumerate}

\end{corollary}
\begin{proof}  (1) $\Rightarrow$ (2): Considering the identity $id_X : X \longrightarrow X$.

(2) $\Rightarrow$ (3): $X^k=P_H(\mathbf{K} (X))=P_H(\mathcal S_c(X)) \cong X$ via a homeomorphism $x\mapsto \overline{\{x\}}$.

(3) $\Rightarrow$ (1): By the adequateness of $\mathbf{K}$, $X^k=P_H(\mathbf{K}(X))$ is a $\mathbf{K}$-space. Since $\mathbf{K}$ is  closed with respect to homeomorphisms and $X\cong X^k$, $X$ is a $\mathbf{K}$-space.
\end{proof}

\begin{corollary}\label{K-retract} Let $\mathbf{K}$ be adequate and closed with respect to homeomorphisms. Then a retract of a $\mathbf{K}$-space is a $\mathbf{K}$-space.
\end{corollary}
\begin{proof} Suppose that $Y$ is a retract of a $\mathbf{K}$-space $X$. Then there are continuous mappings $f : X\longrightarrow Y$ and $g : Y\longrightarrow X$ such that $f\circ g=id_Y$. Let $B\in \mathbf{K}(Y)$, then by Lemma \ref{K-setimage} and Corollary \ref{KKc}, there exists a unique $x_B\in X$ such that $\overline{g(B)}=\overline{\{x_B\}}$. It follows that $B=\overline{f\circ g(B)}=\overline{f(\overline{g(B)})}=\overline{f(\overline{\{x_B\}})}=\overline{\{f(x_B)\}}$. Therefore, $\mathbf{K}(Y)=\mathcal S_c(X)$, and hence $Y$ is a $\mathbf{K}$-space by Corollary \ref{KKc}.
\end{proof}

\begin{theorem}\label{K-reflectionprod}
	For an adequate $\mathbf{K}$ and a finitely family $\{X_i: 1\leq i\leq n\}$ of $T_0$ spaces,  $(\prod\limits_{i=1}^{n}X_i)^k=\prod\limits_{i=1}^{n}X_i^k$ \emph{(}up to homeomorphism\emph{)}.
\end{theorem}

\begin{proof}	
	Let $X=\prod\limits_{i=1}^{n}X_i$. By Corollary \ref{K-closedsetprod}, we can define a mapping $\gamma : P_H(\mathbf{K} (X))  \longrightarrow \prod\limits_{i=1}^{n}P_H(\mathbf{K} (X_i))$ by

\begin{center}
$\forall A\in \mathbf{K} (X)$, $\gamma (A)=(p_1(A), p_2(A), ..., p_n(A))$.
\end{center}

By Lemma \ref{K-setprod} and Corollary \ref{K-closedsetprod}, $\gamma$ is bijective. Now we show that $\gamma$ is a homeomorphism. For any $(U_1, U_2, ..., U_n)\in \mathcal O(X_1)\times \mathcal O(X_2)\times ... \times \mathcal O(X_n)$, by Lemma \ref{K-setprod} and Corollary \ref{K-closedsetprod}, we have

$$\begin{array}{lll}
\gamma^{-1}(\Diamond U_1\times \Diamond U_2\times ... \times\Diamond U_n)&=&\{A\in\mathbf{K}(X): \gamma(A)\in \Diamond U_1\times \Diamond U_2\times ... \times\Diamond U_n\}\\
&=&\{A\in\mathbf{K}(X): p_1(A)\cap U_1\neq\emptyset, p_2(A)\cap U_2\neq\emptyset, ..., p_n(A)\cap U_n\neq\emptyset\}\\
&=&\{A\in\mathbf{K}(X): A\cap U_1\times U_2\times ... \times U_n\neq\emptyset\}\\
&=&\Diamond U_1\times U_2\times ... \times U_n\in \mathcal O(P_H(\mathbf{K} (X)), \mbox{~and~}
\end{array}$$

$$\begin{array}{lll}
\gamma (\Diamond U_1\times U_2\times ... \times U_n)&=&\{\gamma (A): A\in \mathbf{K} (X) \mbox{~and~} A\cap U_1\times U_2\times ... \times U_n\neq\emptyset \}\\
&=&\{\gamma (A): A\in \mathbf{K} (X),\mbox{~and~} p_1(A)\cap U_1\neq\emptyset, p_2(A)\cap U_2\neq\emptyset, ..., p_n(A)\cap U_n\neq\emptyset \}\\
&=&\Diamond U_1\times \Diamond U_2\times ... \times\Diamond U_n\in O(\prod\limits_{i=1}^{n}P_H(\mathbf{K} (X_i))).
\end{array}$$

Therefore, $\gamma : P_H(\mathbf{K} (X))  \longrightarrow \prod\limits_{i=1}^{n}P_H(\mathbf{K} (X_i))$ is a homeomorphism, and hence $X^k$ ($=P_H(\mathbf{K} (X)$) and $\prod\limits_{i=1}^{n}X_i^k$ ($=\prod\limits_{i=1}^{n}P_H(\mathbf{K} (X_i))$ are homeomorphic.
\end{proof}

\begin{theorem}\label{K-prod}
	Suppose that $\mathbf{K}$ is adequate and closed with respect to homeomorphisms. Then for any family $\{X_i:i\in I\}$ of $T_0$ spaces, the following two conditions are equivalent:
	\begin{enumerate}[\rm(1)]
		\item The product space $\prod_{i\in I}X_i$ is a $\mathbf{K}$-space.
		\item For each $i \in I$, $X_i$ is a $\mathbf{K}$-space.
	\end{enumerate}
\end{theorem}
\begin{proof}	
	(1) $\Rightarrow$ (2):  For each $i \in I$, $X_i$ is a retract of $\prod_{i\in I}X_i$. By Corollary \ref{K-retract}, $X_i$ is a $\mathbf{K}$-space.
	
	(2) $\Rightarrow$ (1): Let $X=\prod_{i\in I}X_i$. Suppose $A\in \mathbf{K} (X)$. Then by Corollary \ref{irrcprod}, Corollary \ref{SKIsetrelation} and Lemma \ref{K-setimage}, $A\in \ir_c(X)$ and for each $i \in I$, $p_i(A)\in \mathbf{K} (X_i)$, and consequently, there is a $u_i\in X_i$ such that $p_i(A)=cl_{X_i}\{u_i\}$ by condition (2) and Corollary \ref{KKc}. Let $u=(u_i)_{i\in I}$. Then by Corollary \ref{irrcprod} and \cite[Proposition 2.3.3]{Engelking}),  we have $A=\prod_{i\in I}p_i(A)=\prod_{i\in I}\cl_{u_i}\{u_i\}=\cl_X \{u\}$. Thus $\mathbf{K} (X)=\mathcal{S}_c(X)$. It follows that $X$ is a $\mathbf{K}$-space by Corollary \ref{KKc}.
\end{proof}

\begin{theorem}\label{K-reflection=soberification}
	For an adequate $\mathbf{K}$ and a $T_0$ space $X$, the following conditions are equivalent:
	\begin{enumerate}[\rm (1)]
		\item $X^k$ is the sobrification of $X$, in other words, the $\mathbf{K}$-reflection of $X$ and sobrification of $X$ are the same.
        \item $X^k$ is sober.
		\item $\mathbf{K}(X)=\ir_c(X)$.
	\end{enumerate}
\end{theorem}

\begin{proof}
(1) $\Rightarrow$ (2): Trivial.

(2) $\Rightarrow$ (3): By Corollary \ref{SKIsetrelation}, $\mathbf{K}(X)\subseteq\ir_c(X)$. Now we show that $\ir_c(X)\subseteq\mathbf{K}(X)$. Let $\eta_X^k : X\longrightarrow X^k$ be the canonical topological embedding defined by $\eta_X^k(x)=\overline{\{x\}}$ (see Theorem \ref{K-reflection}). Since the pair $\langle X^s, \eta_{X}^s\rangle$, where $\eta_{X}^s :X\longrightarrow X^s=P_H(\ir_c(X))$, $x\mapsto\overline{\{x\}}$, is the soberification of $X$ and $X^k$ is sober, there exists a unique continuous mapping $(\eta_X^k)^* :X^s\longrightarrow X^k$ such that $(\eta_X^k)^*\circ\eta_X^s=\eta_X^k$, that is, the following diagram commutes.
\begin{equation*}
\xymatrix{
	X \ar[dr]_-{\eta_X^k} \ar[r]^-{\eta_X^s}
	&X^s\ar@{.>}[d]^-{(\eta_X^k)^*}\\
	&X^k}
\end{equation*}	
So for each $A\in\ir_c(X)$, there exists a unique $B\in \mathbf{K} (X)$ such that $\downarrow_{\mathbf{K} (X)} A=\overline{\eta_X^k(A)}=\overline{\{B\}}=\downarrow_{\mathbf{K} (X)}B $. Clearly, we have $B\subseteq A$. On the other hand, for each $a\in A, \overline {\{a\}}\in \downarrow_{\mathbf{K} (X)} A=\downarrow_{\mathbf{K} (X)}B$, and whence $\overline {\{a\}}\subseteq B$. Thus $A\subseteq B$, and consequently, $A=B$. Thus $A\in \mathbf{K}(X)$.

(3) $\Rightarrow$ (1): If $\mathbf{K}(X)=\ir_c(X)$, then $X^k=P_H(\mathbf{K} (X))=P_H(\ir_c(X))=X^s$, with $\eta_X^k=\eta_X^s : X\longrightarrow X^k$, is the sobrification of $X$.

\end{proof}

\begin{proposition}\label{K-reflectioncomp}
	For an adequate $\mathbf{K}$ and a $T_0$ space $X$, $X$ is compact if{}f $X^k$ is compact.
\end{proposition}
\begin{proof} By Corollary \ref{SKIsetrelation}, we have $\mathcal S_c(X)\subseteq \mathbf{K} (X)\subseteq \ir_c (X)$. Suppose that $X$ is compact. For $\{U_i : i\in I\}\subseteq \mathcal O(X)$, if $\mathbf{K} (X)\subseteq \bigcup_{i\in I}\Diamond U_i$, then $X\subseteq \bigcup_{i\in I} U_i$ since $\mathcal S_c(X)\subseteq \mathbf{K} (X)$, and consequently, $X\subseteq \bigcup_{i\in I_0} U_i$ for some $I_0\in I^{(<\omega)}$. It follows that $\mathbf{K} (X)\subseteq \bigcup_{i\in I_0} \Diamond U_i$. Thus $X^k$ is compact. Conversely, if $X^k$ is compact and $\{V_j : j\in J\}$ is a open cover of $X$, then $\mathbf{K} (X)\subseteq \bigcup_{j\in J}\Diamond V_j$. By the compactness of $X^k$, there is a finite subset $J_0\subseteq J$ such that $\mathbf{K} (X)\subseteq \bigcup_{j\in J_0}\Diamond V_j$, and whence $X\subseteq \bigcup_{j\in J_0}V_j$, proving the compactness of $X$.
\end{proof}

For an adequate $\mathbf{K}$, since $\mathcal{S}_c(X)\subseteq \mathbf{K} (X)\subseteq\ir_c(X)$ (see Corollary \ref{SKIsetrelation}), the correspondence $U \leftrightarrow \Diamond U ~(=\Diamond_{\mathbf{K} (X)} U)$ is a lattice isomorphism between $\mathcal O(X)$ and $\mathcal O(X^k)$. Therefore, we have the following result.

\begin{proposition}\label{K-reflectionLHC}
	Let $\mathbf{K}$ be adequate and  $X$  a $T_0$ space. Then
	\begin{enumerate}[\rm (1)]
		\item $X$ is locally hypercompact if{}f $X^k$ is locally hypercompact.
        \item $X$ is a C-space if{}f $X^k$ is a C-space.
        \item $X$ is core compact if{}f $X^k$ is core compact.
	\end{enumerate}
\end{proposition}

\begin{remark}\label{K-reflectionLC} If $\mathbf{K}$ is adequate and $\mathbf{K}\subseteq\mathbf{Top}_w$, then for a $T_0$ space $X$, by Theorem \ref{WFLC=CoreC} and Proposition \ref{K-reflectionLHC}, the following conditions are equivalent:
\begin{enumerate}[\rm (1)]
		\item $X$ is core compact.
        \item $X^k$ is core compact.
        \item $X^k$ is locally compact.
\end{enumerate}
\end{remark}

\begin{remark}\label{corecompnotLC} In \cite{Hofmann-Lawson} (see also \cite[Exercise V-5.25]{redbook}) Hofmann and Lawson given a core compact $T_0$ space $X$ but not locally compact. By Remark \ref{K-reflectionLC} and Theorem \ref{wf-reflection}, $X^s$ and $X^w$ are locally compact. So the local compactness of $X^w$ (or $X^s$) does not imply the local compactness of $X$.
\end{remark}

\begin{definition}\label{Smyth K} $\mathbf{K}$ is said to be a \emph{Smyth category}, if for any $\mathbf{K}$-space $X$, the Smyth power space $P_S(X)$ is a $\mathbf{K}$-space.
\end{definition}

\begin{proposition}\label{WFsober-Smyth} $\mathbf{Sob}$ and $\mathbf{Top}_w$ are Smyth categories.
\end{proposition}
\begin{proof} By \cite[Theorem 3.13]{Klause-Heckmann},  $\mathbf{Sob}$ is a Smyth category. $\mathbf{Top}_w$ is a Smyth category by \cite[Theorem 3]{xuxizhao} or \cite[Theorem 5.3]{xu-shen-xi-zhao-19}.
\end{proof}

\begin{remark}\label{d-not-Smyth} Let $X$ be any $d$-space but not well-filtered (see Example \ref{K-examp1}). Then by \cite[Theorem 5]{xuxizhao}, $P_S(X)$ is not a $d$-space. So $\mathbf{Top}_d$ is not a Smyth category.
\end{remark}

\begin{theorem}\label{SmythK-D}  Let $\mathbf{K}$ be an adequate Smyth category. For a $T_0$ space $X$, if $P_S(X)$ is $\mathbf{K}$-determined, then $X$ is $\mathbf{K}$-determined.
\end{theorem}

\begin{proof} Let $A\in\ir_c(X)$, $Y$ a $\mathbf{K}$-space and $f:X\longrightarrow Y$  a continuous mapping. Then $\overline{\xi_X(A)}=\Diamond A\in\ir_c(P_S(X))$ by Lemma \ref{irrsubspace} and Lemma \ref{irrimage}, and hence $\Diamond A\in \mathbf{K} (P_S(X))$ since $P_S(X)$ is $\mathbf{K}$-determined, where $\xi_X : X \longrightarrow P_S(X)$, $x\mapsto\ua x$. Define a mapping $P_S(f): P_S(X)\longrightarrow P_S(Y)$ by
	$$\forall K\in\mathsf{Q}(X),\ P_S(f)(K)=\ua f(K).$$
	
	{Claim 1:} $P_S(f)\circ\xi_X=\xi_Y\circ f$.
	
	For each $ x\in X$, we have
	$$P_S(f)\circ\xi_X(x)=P_S(f)(\ua x)=\ua f(x)=\xi_Y\circ f(x),$$
	that is, the following diagram commutes.
	\begin{equation*}
	\xymatrix{
		X \ar[d]_-{f} \ar[r]^-{\xi_X} &P_S(X)\ar[d]^-{P_S(f)}\\
		Y \ar[r]^-{\xi_Y} &P_S(Y)	}
	\end{equation*}

{Claim 2:} $P_S(f): P_S(X)\longrightarrow P_S(Y)$ is continuous.

Let $V\in\mathcal O(Y)$. We have
$$\begin{array}{lll}
P_S(f)^{-1}(\Box V)&=&\{K\in\mathsf{Q}(X): P_S(f)(K)=\ua f(K)\subseteq V\}\\
&=&\{K\in\mathsf{Q}(X): K\subseteq f^{-1}(V)\}\\
&=&\Box f^{-1}(V),
\end{array}$$
which is open in $P_S(X)$. This implies that $P_S(f)$ is continuous.

Since $\mathbf{K}$ is a Smyth category, $P_S(Y)$ is a $\mathbf{K}$-space. By the continuity of  $P_S(f)$ and $\Diamond A\in \mathbf{K}(P_S(X))$,
there exists a unique $Q\in \mathsf{Q}(Y)$ such that $\overline{P_S(f)(\Diamond A)}=\overline{\{Q\}}$.

 {Claim 3:} $Q$ is supercompact.

 Let $\{U_j:j\in J\}\subseteq\mathcal O(X)$ with $Q\subseteq \bigcup_{j\in J}U_j$, i.e., $Q\in\Box \bigcup_{j\in J}U_j$. Note that $\overline{P_S(f)(\Diamond A)}=\overline{\{\ua f(a):a\in A\}}$, thus $\{\ua f(a):a\in A\}\cap \Box \bigcup_{j\in J}U_j\neq\emptyset$. Then there exists $a_0\in A$ and $j_0\in J$ such that
 $Q\subseteq \ua f(a_0)\subseteq U_{j_0}$.

 Hence, by \cite[Fact 2.2]{Klause-Heckmann}, there exists $y_Q\in Y$ such that $Q=\ua y_Q$.

 {Claim 4:} $\overline{f(A)}=\overline{\{y_Q\}}$.

 Note that $\overline{\{\ua f(a):a\in A\}}=\overline{\{\ua y_Q\}}$. Thus for each $y\in f(A)$, $\ua y\in \overline{\{\ua y_Q\}}$, showing that $\ua y_Q\subseteq \ua y$, i.e.,  $y\in\overline{\{y_Q\}}$. This implies that $f(A)\subseteq \overline{\{y_Q\}}$.
 In addition, since $\ua y_Q\in\overline{\{\ua f(a):a\in A\}}=\Diamond \overline{f(A)}$, $\ua y_Q\cap\overline{f(A)}\neq\emptyset$. This implies that $y_Q\in \overline{f(A)}$m, and whence $\overline{f(A)}=\overline{\{y_Q\}}$. Thus $A\in \mathbf{K}(X)$. Therefore, by Corollary \ref{SKIsetrelation}, $\mathbf{K}(X)=\ir_c(X)$, proving that $X$ is $\mathbf{K}$-determined.
 \end{proof}

\section{Applications}

This section is devoted to giving some applications of the results of Section 4 to $\mathbf{Sob}$, $\mathbf{Top}_d$, $\mathbf{Top}_w$ and the Keimel-Lawson category.

First, we consider the case of $\mathbf{K}=\mathbf{Sob}$. For a $T_0$ space $X$, by Lemma \ref{sobd=irr}, $\mathbf{Sob}(X)=\ir_c(X)$. It is well-known that $P_H(\ir_c(X))$ is sober (see, e.g., \cite{redbook, Jean-2013}). In fact, for any $\mathcal A\in \ir_c(P_H(\ir_c(X)))$, let $A=\overline{\cup \mathcal A}$. Then $A\in \ir_c(X)$ and $\mathcal A=\overline{\{A\}}$ in $P_H(\ir_c(X))$. Thus $P_H(\ir_c(X))$ is sober. Therefore, by Proposition \ref{WFsober-Smyth}, we get the following well-known result.

\begin{proposition}\label{sobrification}  $\mathbf{Sob}$ is an adequate Smyth category. Therefore, for any $T_0$ space $X$, $X^s=P_H(\ir_c(X))$ with the canonical mapping $\eta_{X}: X\longrightarrow X^s$ is the \emph{sobrification} of $X$.
\end{proposition}

It follows from Proposition \ref{sobrification} that $\mathbf{Sob}$ is reflective in $\mathbf{Top}_0$ (cf. \cite{redbook}).

 \begin{proposition}\label{K-reflectionprod}\emph{(\cite{Hoffmann1, Hoffmann2})}
	For a family $\{X_i: i\in I\}$ of $T_0$ spaces, $(\prod_{i\in I}X_i)^s=\prod_{i\in I}X_i^s$ \emph{(}up to homeomorphism\emph{)}.
\end{proposition}

\begin{proof}	
	Let $X=\prod_{i\in I}X_i$. By Lemma \ref{prodirr} and Corollary \ref{irrcprod}, we can define a bijective mapping $\beta : P_H(\ir_c(X))  \longrightarrow \prod_{i\in I}P_H(\ir_c(X_i))$ by

\begin{center}
$\forall A\in \ir_c(X)$, $\beta (A)=(p_i(A))_{i\in I}$.
\end{center}

  Now we show that $\beta$ is a homeomorphism. Let $q_i : \prod_{i\in I}P_H(\ir_c(X_i)) \longrightarrow P_H(\ir_c(X_i))$ be the $i$th projection ($i\in I$).  For any $J\in I^{(<\omega)}$ and $(U_i)_{i\in J}\in \prod_{i\in J}\mathcal O(X_i)$, by Lemma \ref{prodirr} and Corollary \ref{irrcprod}, we have

$$\begin{array}{lll}
\beta^{-1}(\bigcap_{i\in J}q_i^{-1}(\Diamond U_i))&=&\{A\in\ir_c(X): \beta(A)\in \bigcap_{i\in J}q_i^{-1}(\Diamond U_i)\}\\
&=&\{A\in\ir_c(X): p_i(A)\cap U_i\neq\emptyset \mbox{ for each } i\in J \}\\
&=&\{A\in\ir_c(X): A\cap \bigcap_{i\in J}p_i^{-1}( U_i)\neq\emptyset\}\\
&=&\Diamond \bigcap_{i\in J}p_i^{-1} (U_i)\in \mathcal O(P_H(\ir_c (X)), \mbox{~and~}
\end{array}$$

$$\begin{array}{lll}
\beta (\Diamond \bigcap_{i\in J}p_i^{-1} (U_i))&=&\{\beta (A): A\in \ir_c (X) \mbox{~and~} A\cap \bigcap_{i\in J}p_i^{-1}( U_i)\neq\emptyset \}\\
&=&\bigcap_{i\in J}p_i^{-1}(\Diamond U_i)\in O(\prod\limits_{i\in I}P_H(\ir_c (X_i))).
\end{array}$$

Therefore, $\beta : P_H(\ir_c (X))  \longrightarrow \prod_{i\in I}P_H(\ir_c (X_i))$ is a homeomorphism, and hence $X^s$ ($=P_H(\ir_c (X)$) and $\prod\limits_{i\in I}X_i^s$ ($=\prod_{i\in I}P_H(\ir_c (X_i))$ are homeomorphic.
\end{proof}

By Theorem \ref{K-prod} and Proposition \ref{sobrification}, we get the following well-known result (see, e.g. \cite[Exercise O-5.16]{redbook}).

\begin{corollary}\label{sober-prod}
	For a family $\{X_i:i\in I\}$ of $T_0$ spaces, the following two conditions are equivalent:
	\begin{enumerate}[\rm(1)]
		\item The product space $\prod_{i\in I}X_i$ is sober.
		\item For each $i \in I$, $X_i$ is sober.
	\end{enumerate}
\end{corollary}

Second, we consider the case of $\mathbf{K}=\mathbf{Top}_d$.

\begin{theorem}\label{d-reflection}  $\mathbf{Top}_d$ is adequate. Therefore, for any $T_0$ space $X$, $X^d=P_H(d(X))$ with the canonical mapping $\eta_{X}: X\longrightarrow X^d$ is the $d$-reflection of $X$.
\end{theorem}

\begin{proof}
	Suppose that $X$ be a $T_0$ space. We show that $P_H(d(X))$ is a $d$-space. Since $X$ is $T_0$, one can directly deduce that $P_H(d(X))$ is $T_0$. Let $\{
	A_d:d\in D\}\subseteq d(X)$ be a directed family. Let $A=\overline{\bigcup_{d\in D}A_d}$. We check that $A\in d(X)$. For any continuous mapping $f : X \longrightarrow Y$ to a $d$-space $Y$ and $d\in D$, by $A_d\in d(X)$, there is a $y_d\in Y$ such that $\overline{f(A_d)}=\overline{\{y_d\}}$. Since $\{
	A_d:d\in D\}\subseteq d(X)$ is directed, $\{y_d : d\in D\}\in \mathcal D(Y)$. By Proposition \ref{d-spacecharac1}, there is a $y\in Y$ such that $\overline{\{y_d : d\}}=\overline{\{y\}}$.  Therefore, we have
 $$\begin{array}{lll}
           \overline{f(A)} & =\overline{f(\overline{\bigcup_{d\in D}A_d})}\\
	                       & =\overline{f(\bigcup_{d\in D}A_d)}\\
	                       & =\overline{\bigcup_{d\in D}f(A_d)}\\
	                       & =\overline{\bigcup_{d\in D}\overline{f(A_d)}}\\
                           & =\overline{\bigcup_{d\in D}\overline{\{y_d\}}}\\
                           & =\overline{\{y_d : d\in D\}}\\
                           & =\overline{\{y\}}.\\
	\end{array}$$
Thus $A\in d(X)$. Clearly, $\overline{\{
	A_d:d\in D\}}=\overline{\{A\}}$ in $P_H(d(X))$. By Proposition \ref{d-spacecharac1} again, $P_H(d(X))$ is a $d$-space.

By Lemma \ref{K-lemmafstar}, the pair $\langle X^d=P_H(d(X)), \eta_X^d\rangle$, where $\eta_X^d :X\longrightarrow X^d$, $x\mapsto\overline{\{x\}}$, is the $d$-reflection of $X$.
\end{proof}

\begin{corollary}\label{d-reflective}\emph{(\cite{Ershov_1999, Keimel-Lawson, Wyler})}  $\mathbf{Top}_d$ is reflective in $\mathbf{Top}_0$.
\end{corollary}

From Theorem \ref{K-prod} and Theorem \ref{d-reflection} we deduce the following known result.

\begin{corollary}\label{d-prod}
	For a family $\{X_i:i\in I\}$ of $T_0$ spaces, the following two conditions are equivalent:
	\begin{enumerate}[\rm(1)]
		\item The product space $\prod_{i\in I}X_i$ is a $d$-space.
		\item For each $i \in I$, $X_i$ is a $d$-space.
	\end{enumerate}
\end{corollary}

By Theorem \ref{K-reflectionprod} and Theorem \ref{d-reflection}, we get the following two results, which were  proved by Keimel and Lawson using $d$-closures in \cite{Keimel-Lawson}.

\begin{corollary}\label{d-reflectionprod} \emph{(\cite{Keimel-Lawson})}
	For a finitely family $\{X_i: 1\leq i\leq n\}$ of $T_0$ spaces,  $(\prod\limits_{i=1}^{n}X_i)^d=\prod\limits_{i=1}^{n}X_i^d$ \emph{(}up to homeomorphism\emph{)}.
\end{corollary}

\begin{corollary}\label{d-reflectionScott} \emph{(\cite{Keimel-Lawson})}
	For a $T_0$ space $X$, if $\xi_X : X \longrightarrow \Sigma~\!\! d(X)$, $\xi_X(x)=\overline{\{x\}}$, is continuous, then the $d$-reflection $X^d$ of $X$ is the Scott space $\Sigma~\!\! d(X)$.
\end{corollary}

\begin{proof} By Theorem \ref{d-reflection}, $X^d=P_H(d(X))$ with $\eta_X : X \longrightarrow X^d$ is the $d$-reflection of $X$, and consequently, $d(X)$ (with respect to the specialization order or, equivalently, the order of set inclusion) is a dcpo, $\Sigma~\!\! d(X)$ is a $d$-space, and $\mathcal O(X^d)\subseteq \sigma (d(X))$. Since $\xi_X : X \longrightarrow \Sigma~\!\! d(X)$ is continuous, there is a unique continuous mapping $(\xi_X)^d : P_H(d(X))\longrightarrow \Sigma~\!\! d(X)$ such that $(\xi_X)^d\circ\eta_X=\xi_X$, that is, the following diagram commutes.
\begin{equation*}
\xymatrix{
	X \ar[dr]_-{\xi_X} \ar[r]^-{\eta_X}
	&P_H(d(X))\ar@{.>}[d]^-{(\xi_X)^d}\\
	&\Sigma~\!\! d(X)}
\end{equation*}	

For each $A\in d(X)$, by Lemma \ref{lemmaclosure} and the proof of Lemma \ref{K-lemmafstar}, there exists a unique $B\in d(X)$ such that $\da_{d(X)}A=\Box A=\overline{\xi_X(A)}=\overline{\{B\}}=\da_{d(X)}B$. Therefore, $A=B$ (note that $\mathcal{S}_c(X)\subseteq d(X)$), and hence $(\xi_X)^d(A)=A$. It follows that $\sigma(d(X))\subseteq \mathcal O(X^d)$. Thus $\mathcal O(X^d)=\sigma (d(X))$.

\end{proof}

The following result shows that the $d$-reflection space of Scott space of a poset $P$ is the Scott space $\Sigma~\!\! d(\Sigma~\!\! P))$ of the dcpo $d(\Sigma~\!\! P)$.

\begin{corollary}\label{d-reflectionScott}
	For any poset $P$, $d(\Sigma~\!\! P)$ is a dcpo and the $d$-reflection  $(\Sigma~\!\! P)^d$ of $\Sigma~\!\! P$ is the Scott space $\Sigma~\!\! d(\Sigma~\!\! P))$ with the canonical mapping $\eta_{P}: \Sigma~\!\! P\longrightarrow \Sigma~\!\! d(X)$, given by $\eta_{P}(x)=cl_{\sigma(P)}\{x\}$ for each $x\in X$.
\end{corollary}
\begin{proof} By Corollary \ref{continuous-ScottCONT}, Theorem \ref{d-reflection} and Corollary \ref{d-reflectionScott}.
\end{proof}

\begin{corollary}\label{posets-reflection}\emph{(\cite{ZhaoFan})}
	$\mathbf{DCPO}$ is reflective in $\mathbf{Poset}_s$.
\end{corollary}

\begin{definition}\label{DCPO-comp}(\cite{ZhaoFan})
	A $\mathbf{DCPO}$-\emph{completion} of a poset $P$, $\mathbf{D}$-completion of $P$ for short, is a pair $\langle \widetilde{P}, \eta\rangle$ consisting of a dcpo  $\widetilde{P}$ and a Scott continuous mapping $\eta :P\longrightarrow \widetilde{P}$, such that for any Scott continuous mapping $f: P\longrightarrow Q$ to a dcpo $Q$, there exists a unique Scott continuous mapping $\widetilde{f} : \widetilde{P}\longrightarrow Q$ such that $\widetilde{f}\circ\eta=f$, that is, the following diagram commutes.\\
\begin{equation*}
	\xymatrix{
		P \ar[dr]_-{f} \ar[r]^-{\eta}
		&\widetilde{P}\ar@{.>}[d]^-{\widetilde{f}}\\
		&Q}
	\end{equation*}

\end{definition}

$\mathbf{D}$-completions, if they exist, are unique up to isomorphism. We shall use $\mathbf{D}(P)$ to denote the $\mathbf{D}$-completion of $P$ if it exists.

In \cite{ZhaoFan}, using the $D$-topologies defined in \cite{ZhaoFan}, Zhao and Fan proved that for any poset $P$, the $\mathbf{D}$-completion of $P$ exists. As Keimel and Lawson  pointed out in \cite{Keimel-Lawson} that the $\mathbf{D}$-completion of a poset $P$ is a special case of the $d$-reflection of a certain $T_0$ space. More precisely, the $d$-reflection of Scott space $\Sigma~\!\! P$.

\begin{proposition}\label{DCPO-completion} For a poset $P$, $\mathbf{D}(P)=d(\Sigma~\!\! P)$ with the canonical mapping $\eta_{P}: P\longrightarrow \mathbf{D}(P)$, $\eta_P(x)=cl_{\sigma(P)}\{x\}$, is the $\mathbf{D}$-completion of $P$.
\end{proposition}
\begin{proof} By Theorem \ref{d-reflection},  $(\Sigma~\!\! P)^d=P_H(d(\Sigma~\!\! P))$ with the canonical mapping $\eta_{P}: \Sigma~\!\! P\longrightarrow (\Sigma~\!\! P)^d$, $\eta_P(x)=cl_{\sigma(P)}\{x\}$, is the $d$-reflection of $\Sigma~\!\! P$. By Lemma \ref{Scott-cont1} and Corollary \ref{d-reflectionScott}, $\mathbf{D}(P)=d(\Sigma~\!\! P)$ with the canonical mapping $\eta_{P}: P\longrightarrow \mathbf{D}(P)$ is the $\mathbf{D}$-completion of $P$.
\end{proof}

Now we consider the case of $\mathbf{K}=\mathbf{Top}_w$.

\begin{lemma}\label{lemmaBoxK-set}
	Let $X$ be a $T_0$ space and $C\in\mathcal C(X)$. Then the following conditions are equivalent:
	\begin{enumerate}[\rm (1)]
		\item $C\in \mathbf{K}(X)$.
		\item $\Box C\in \mathbf{K}(P_H(\mathbf{K}(X)))$.
	\end{enumerate}
\end{lemma}
\begin{proof}
(1) $\Rightarrow$ (2): By Propositions \ref{K-setimage}, Lemma \ref{lemmaeta} and Lemma \ref{lemmaclosure}.

(2) $\Rightarrow$ (1).  Let $Y$ be a $\mathbf{K}$-space and $f:X\longrightarrow Y$  a continuous mapping. By Lemma \ref{K-lemmafstar}, there exists a continuous mapping $f^* :P_H(\mathbf{K}(X))\longrightarrow Y$ such that $f^*\circ\eta_X=f$.
Since $\Box C=\overline{\eta_X(C)}$ is a $\mathbf{K}$-set and $f^*$ is continuous,  there exists a unique $y_C\in Y$ such that
$\overline{f^*\left(\overline{\eta_X(C)}\right)}=\overline{\{y_C\}}$. Furthermore, we have
$$\overline{\{y_C\}}=\overline{f^*\left(\overline{\eta_X(C)}\right)}=\overline{f^*(\eta_X(C))}=\overline{f(C)}.$$
So $C$ is a $\mathbf{K}$-set.
\end{proof}

\begin{theorem}\label{wf-reflection}  $\mathbf{Top}_w$ is adequate. Therefore, for any $T_0$ space $X$, $X^w=P_H(\mathsf{WF}(X))$ with the canonical mapping $\eta_{X}: X\longrightarrow X^w$ is the well-filtered reflection of $X$.
\end{theorem}
\begin{proof}
	Suppose that $X$ be a $T_0$ space. We show that $P_H(\mathsf{WF}(X))$ is well-filtered. Since $X$ is $T_0$, one can directly check that $P_H(\mathsf{WF}(X))$ is $T_0$. Let $\{
	\mathcal K_i:i\in I\}\subseteq\mathsf{Q}(P_H(\mathsf{WF}(X)))$ be a filtered family and $U\in\mathcal O(X)$ such that
	$\bigcap_{i\in I}\mathcal K_i\subseteq \Diamond U$. We need to show $\mathcal K_i\subseteq\Diamond U$ for some $i\in I$.
	Assume, on the contrary, $\mathcal K_i\nsubseteq \Diamond U$, i.e., $\mathcal K_i\cap\Box(X\setminus U)\neq\emptyset$,  for any $i\in I$.
	
	Let $\mathcal A=\{C\in\mathcal C(X): C\subseteq X\setminus U \text{ and } \mathcal K_i\cap \Box C\neq\emptyset \text{ for all }  i\in I\}$. Then we have the following two facts.
	
	{(a1)} $\mathcal A\neq\emptyset$ because $X\setminus U\in\mathcal A$.
	
	{(a2)} For any filtered family $\mathcal F\subseteq\mathcal A$, $\bigcap\mathcal F\in\mathcal A$.
	
	Let $F=\bigcap\mathcal F$. Then $F\in \mathcal C(X)$ and $F\subseteq X\setminus U$. Assume, on the contrary, $F\notin\mathcal A$. Then there exists $i_0\in I$ such that $\mathcal K_{i_0}\cap \Box F=\emptyset$. Note that $\Box F=\bigcap_{C\in\mathcal F}\Box C$, implying that $\mathcal K_{i_0}\subseteq\bigcup_{C\in\mathcal F}\Diamond (X\setminus C)$ and $\{\Diamond (X\setminus C) : C\in\mathcal F\}$ is a directed family since $\mathcal F$ is filtered. Then there is $C_0\in\mathcal F$ such that $\mathcal K_{i_0}\subseteq \Diamond (X\setminus C_0)$, i.e., $\mathcal K_{I_0}\cap\Box C_0=\emptyset$,  contradicting $C_0\in\mathcal A$. Hence $F\in\mathcal A$.
	
	By Zorn's Lemma, there exists a minimal element $C_m$ in $\mathcal A$ such that $\Box C_m$ intersects all members of $\mathcal K$. Clearly,  $\Box C_m$ is also a minimal closure set that  intersects all members of $\mathcal K$, hence is a Rudin set in $P_H(\mathbf{K}(X))$. By Proposition \ref{DRWIsetrelation} and  Lemma \ref{lemmaBoxK-set}, $C_m\in \mathsf{WF}(X)$. So $C_m\in \Box C_m\cap\bigcap\mathcal K\neq\emptyset$. It follows that $\bigcap\mathcal K\nsubseteq \Diamond(X\setminus C_m)\supseteq \Diamond U$, which implies that $\bigcap\mathcal K\nsubseteq \Diamond U$, a contradiction.

By Lemma \ref{K-lemmafstar}, the pair $\langle X^w=P_H(\mathsf{WF}(X)), \eta_X^w\rangle$, where $\eta_X^w :X\longrightarrow X^w$, $x\mapsto\overline{\{x\}}$, is the well-filtered reflection of $X$.

\end{proof}

\begin{corollary}\label{wf-reflective}\emph{(\cite{Shenchon, wu-xi-xu-zhao-19, xu-shen-xi-zhao-19})}  $\mathbf{Top}_w$ is reflective in $\mathbf{Top}_0$.
\end{corollary}

By Theorem \ref{K-prod} and Theorem \ref{wf-reflection}, we have the following result.

\begin{corollary}\label{WF-prod}\emph{(\cite{Shenchon, xu-shen-xi-zhao-19})}
	For a family $\{X_i:i\in I\}$ of $T_0$ spaces, the following two conditions are equivalent:
	\begin{enumerate}[\rm(1)]
		\item The product space $\prod_{i\in I}X_i$ is well-filtered.
		\item For each $i \in I$, $X_i$ is a well-filtered.
	\end{enumerate}
\end{corollary}

Finally, we consider the case that $\mathbf{K}$ is a Keimel-Lawson category.

\begin{theorem}\label{KL-reflection} Let $\mathbf{K}$ be a  Keimel-Lawson category. Then $\mathbf{K}$ is adequate. Therefore, for any $T_0$ space $X$, $X^k=P_H(\mathbf{K}(X))$ with the canonical mapping $\eta_{X}^k: X\longrightarrow X^k$ is the $\mathbf{K}$-reflection of $X$.
\end{theorem}

\begin{proof} For any $\mathbf{K}$-space $Y$ and continuous mapping $f : X \longrightarrow Y$, let $j_X^k : X^k \longrightarrow X^s$ be the inclusion mapping (note that $\mathbf{K}(X)\subseteq \ir_c(X)$). By Lemma \ref{K-lemmafstar} and Proposition \ref{sobrification}, there is a unique continuous mappings
$f^k : X^k \longrightarrow Y$ such that $f^k\circ \eta_X^k=f$, and a unique continuous mapping $f^* : X^s \longrightarrow Y^s$ such that $f^*\circ \j_X^k=\eta_Y\circ f^k$, that is, the following diagram commutes.
		\begin{equation*}
	\xymatrix{
		X \ar[d]_-{f} \ar[r]^-{\eta_X^k} &X^k\ar[d]^-{f^k} \ar[r]^-{j_X^k} &X^s\ar[d]^-{f^*} \\
		Y \ar[r]^-{id_Y} &Y \ar[r]^-{\eta_Y^s} &Y^s
	}
	\end{equation*}
For each $A\in \mathbf{K} (X)$, $\overline{\{f^k(A)\}}=\overline{\{f^*(A)\}}=\overline{f(A)}$. Let $[X\rightarrow\mathbf{K}]=\{f : X \longrightarrow Y~ | ~Y$ is a $\mathbf{K}$-space and $f$ is continuous$\}$. We have that $\mathbf{K}(X)=\bigcap\limits_{f\in [X\rightarrow\mathbf{K}]}(f^*)^{-1}(\{\overline{\{y\}} : y\in Y\})$, and whence $P_H(\mathbf{K}(X))=\bigcap\limits_{f\in [X\rightarrow\mathbf{K}]}P_H((f^*)^{-1}(P_H(\mathcal S_c(Y))))$ (see Remark \ref{PH-subtopology}). $P_H(\mathcal S_c(Y))$ is homeomorphic to $Y$, and hence it is a $\mathbf{K}$-space by ($\mathrm{K}_1$). For each $f\in [X\rightarrow\mathbf{K}]$, by ($\mathrm{K}_1$) and ($\mathrm{K}_4$), $P_H((f^*)^{-1}(P_H(\mathcal S_c(Y))))$ (as a subspace of $X^s$) is a $\mathbf{K}$-space. Finally, by ($\mathrm{K}_3$), $P_H(\mathbf{K}(X))$ is a $\mathbf{K}$-space. Thus $\mathbf{K}$ is adequate. By Lemma \ref{K-lemmafstar}, the pair $\langle X^k=P_H(\mathbf{K}(X)), \eta_X^w\rangle$, where $\eta_X^k :X\longrightarrow X^k$, $x\mapsto\overline{\{x\}}$, is  the $\mathbf{K}$-reflection of $X$.
\end{proof}

\begin{corollary}\label{K-reflective} \emph{(\cite{Keimel-Lawson})} Every Keimel-Lawson category $\mathbf{K}$ is reflective in $\mathbf{Top}_0$.
\end{corollary}

From Theorem \ref{K-reflectionprod} and Theorem \ref{KL-reflection} we deduce the following corollary.

\begin{corollary}\label{d-reflectionprod}
	Let $\mathbf{K}$ be a Keimel-Lawson category. For a finitely family $\{X_i: 1\leq i\leq n\}$ of $T_0$ spaces,  $(\prod\limits_{i=1}^{n}X_i)^k=\prod\limits_{i=1}^{n}X_i^k$ \emph{(}up to homeomorphism\emph{)}.
\end{corollary}

By Theorem \ref{K-prod} and Theorem \ref{KL-reflection}, we get the following result.

\begin{corollary}\label{KL-prod}
	Let $\mathbf{K}$ be a Keimel-Lawson category. Then for a family $\{X_i:i\in I\}$ of $T_0$ spaces, the following two conditions are equivalent:
	\begin{enumerate}[\rm(1)]
		\item The product space $\prod_{i\in I}X_i$ is a $\mathbf{K}$-space.
		\item For each $i \in I$, $X_i$ is a $\mathbf{K}$-space.
	\end{enumerate}
\end{corollary}

\section{Conclusion}

In this paper, we provided a direct approach to $\mathbf{K}$-reflections of $T_0$ spaces. For a full subcategory $\mathbf{K}$ of $\mathbf{Top}_0$ containing $\mathbf{Sob}$ and a $T_0$ space $X$, it was proved that if $P_H(\mathbf{K}(X))$ is a $\mathbf{K}$-space,  then the pair $\langle X=P_H(\mathbf{K}(X)), \eta_X\rangle$, where $\eta_X :X\longrightarrow X^k$, $x\mapsto\overline{\{x\}}$, is the $\mathbf{K}$-reflection of $X$. Therefore,  every adequate $\mathbf{K}$ is reflective in $\mathbf{Top}_0$. It was shown that $\mathbf{Sob}$, $\mathbf{Top}_d$, $\mathbf{Top}_w$ and the Keimel and Lawson's category are all adequate, and hence they are all reflective in $\mathbf{Top}_0$. Some major properties of  $\mathbf{K}$-spaces and $\mathbf{K}$-reflections of $T_0$ spaces were investigated. In particular, it was proved that if $\mathbf{K}$ is adequate, then the $\mathbf{K}$-reflection preserves finite products of $T_0$ spaces. Our study also leads to a number of problems, whose answering will deepen our understanding of the related spaces and their categorical structures.

In \cite{ZhaoHo}, Zhao and Ho introduced a weak notion of sobriety: $k$-bounded sobriety. Recently, Ern\'e \cite{E_20182} replaced joins by cuts,
and introduced three kinds of non-sober spaces: cut spaces, weakly sober spaces, and quasisober spaces. In
a forthcoming article we will show that some of the categories of $k$-bounded sober spaces, cut spaces, weakly sober spaces, and quasisober spaces are not adequate and they are really not reflective in $\mathbf{Top}_0$.

We now close our paper with the following questions about $\mathbf{K}$-reflections of $T_0$ spaces, where $\mathbf{K}$ is an adequate full subcategory of $\mathbf{Top}_0$ containing $\mathbf{Sob}$.

\begin{question}\label{KL-reflection-infinite-prod} Does the $\mathbf{K}$-reflection (especially, for a Keimel-Lawson category $\mathbf{K}$) preserve arbitrary products of $T_0$ spaces? Or equivalently, does $(\prod\limits_{i\in I}X_i)^k=\prod\limits_{i\in I}X_i^k$ (up to homeomorphism) hold for any family $\{X_i : i\in I\}$ of $T_0$ spaces?
\end{question}

\begin{question}\label{K-reflection-infinite-prod} Does the $d$-reflection preserve arbitrary products of $T_0$ spaces?
\end{question}

\begin{question}\label{d-reflection-infinite-prod} Does the well-filtered reflection preserve arbitrary products of $T_0$ spaces?
\end{question}

\begin{question}\label{K-infiniteset-prodquestion} Let $X=\prod_{i\in I}X_i$ be the product space of a family $\{X_i: i\in I\}$ of $T_0$ spaces. If each $A_i\subseteq X_i~(i\in I)$ is a $\mathbf{K}$-set, must the product set $\prod_{i\in I}A_i$ be a $\mathbf{K}$-set of $X$?
\end{question}

\begin{question}\label{WDinfinite-prodquestion} Is the product space of an arbitrary  collection of $\mathbf{K}$-determined spaces $\mathbf{K}$-determined?
\end{question}

\begin{question}\label{SmythWD} Let $\mathbf{K}$ be a Smyth category. Is the Smyth power space $P_S(X)$ of a $\mathbf{K}$-determined $T_0$ space $X$ again $\mathbf{K}$-determined?
\end{question}

\end{document}